\documentclass[12pt]{amsart}
\usepackage{amsfonts}
\usepackage{amsmath}
\usepackage{amssymb}
\usepackage[mathcal]{eucal}
\usepackage{amscd}
\input xy
\xyoption{all}

\newcommand{\Ucal}{\mathcal{U}}

\newcommand{\Mf}{\mathfrak{M}}

\newcommand{\R}{\mathbb{R}}

\newcommand{\N}{\mathbb{N}}

\newcommand{\XtoYGa}{$(X,\Ga)\overset\pi\to(Y,\Ga)$\ }

\newcommand{\al}{\alpha}
\newcommand{\Ga}{\Gamma}
\newcommand{\ga}{\gamma}
\newcommand{\del}{\delta}
\newcommand{\Del}{\Delta}
\newcommand{\ep}{\epsilon}
\newcommand{\sig}{\sigma}

\newcommand{\la}{\lambda}
\newcommand{\La}{\Lambda}
\newcommand{\tet}{\theta}

\newcommand{\ka}{\kappa}

\newcommand{\XGa}{$(X,\Ga)\ $}
\newcommand{\YGa}{$(Y,\Ga)\ $}

\newcommand{\OGa}{\mathcal{O}_{\Ga}}
\newcommand{\OCGa}{ {\ol{\mathcal{O}}}_{\Ga} }

\newcommand{\br}{\vspace{4 mm}}

\newcommand{\rest}{\upharpoonright}

\newcommand{\ol}{\overline}
\newcommand{\ch}{\mathbf{1}}

\newcommand{\inte}{\rm{inter\,}}
\newcommand{\cls}{\rm{cls\,}}

\newcommand{\Aut}{\rm{Aut\,}}

\newcommand{\gr}{\rm{graph\,}}

\newcommand{\diam}{\rm{diam\,}}
\newcommand{\supp}{\rm{supp\,}}

\newcommand{\card}{\rm{card\,}}

\newcommand{\Homeo}{\rm{Homeo\,}}

\swapnumbers
\theoremstyle{plain}
\newtheorem{thm}{Theorem}[section]

\newtheorem{lem}[thm]{Lemma}
\newtheorem{prop}[thm]{Proposition}

\theoremstyle{definition}

\newtheorem{rmk}[thm]{Remark}

\numberwithin{equation}{section}

\begin{document}
\title
[Tame minimal dynamical systems]{The structure of tame minimal
dynamical systems}

\author[Eli Glasner]{Eli Glasner}
\address{Department of Mathematics,
Tel-Aviv University, Tel Aviv, Israel}
\email{glasner@math.tau.ac.il}

\begin{abstract}
A dynamical version of the Bourgain-Fremlin-Talagrand
dichotomy shows that the enveloping semigroup of a
dynamical system is either very large and contains
a topological copy of $\beta \N$, or it is a
``tame" topological space whose topology is determined by
the convergence of sequences. In the latter case
the dynamical system is called tame.
We use the structure theory of minimal dynamical systems
to show that, when the acting group is Abelian,
a tame metric minimal dynamical system (i) is
almost automorphic (i.e. it is an almost 1-1 extension
of an equicontinuous system), and (ii) admits a unique
invariant probability measure such that the corresponding
measure preserving system is measure theoretically isomorphic
to the Haar measure system on the maximal equicontinuous factor.
\end{abstract}

\thanks{{\em 2000 Mathematical Subject Classification
54H20}}


\maketitle

\tableofcontents

\section*{Introduction}
In this work a dynamical system is a pair $(X,\Ga)$,
where $X$ is a compact Hausdorff space and $\Ga$ an
abstract group acting as a group of homeomorphisms of the space $X$.
That is we are given a homomorphism (not necessarily an
isomorphism) of $\Ga$ into $\Homeo(X)$.
For $\ga\in \Ga$ and $x\in X$ we write $\ga x$ for the image
of $x$ under the homeomorphism which corresponds to $\ga$.
We will often abuse this notation and consider $\ga$ as
a homeomorphism of $X$.

The {\em enveloping semigroup} $E(X,\Ga)$ of the dynamical system \XGa
is defined as the closure of image of $\Ga$ in the product space $X^X$.
It is not hard to check that, under composition of maps,
$E(X,\Ga)$ is a compact {\em right topological semigroup} , i.e.
for each $q\in E(X,\Ga)$ the map $R_q:p \mapsto pq$ is continuous.
In fact the canonical map of $\Ga$ into $E(X,\Ga)$ is a
{\em right topological semigroup compactification} of $\Ga$; i.e.
it has a dense range and for each $\ga \in \Ga$
multiplication on the left $L_\ga: p \mapsto \ga p$
is continuous on $E(X,\Ga)$. This left multiplication by elements of $\Ga$
makes $(E(X,\Ga),\Ga)$ a dynamical system.

The enveloping semigroup was introduced by Robert Ellis in 1960 and
became an indispensable tool in abstract topological
dynamics. However explicit computations of enveloping
semigroups are quite rare.
One reason for this is that often $E(X,G)$ is non-metrizable.

Following an idea of A. K\"{o}ller, \cite{Ko},
Glasner and Megrelishvily proved the
following dynamical version of the
Bourgain-Fremlin-Talagrand dichotomy theorem, \cite{GM}.

\begin{thm}[A dynamical BFT dichotomy]
Let $(X,G)$ be a metric dynamical system and let $E(X)$
be its enveloping semigroup. We have the following dichotomy.
Either
\begin{enumerate}
\item
$E(X)$ is separable Rosenthal compact, hence with cardinality
${\card}\ {E(X)} \leq 2^{\aleph_0}$; or
\item
the compact space $E$ contains a homeomorphic
copy of $\beta\mathbb{N}$, hence ${\card}\ {E(X)} = 2^{2^{\aleph_0}}$.
\end{enumerate}
\end{thm}

A dynamical system is called {\em tame} if the first alternative occurs,
i.e. $E(X)$ is Rosenthal compact.
Recently dynamical characterizations of both tame dynamical systems
and dynamical systems whose enveloping groups are metrizable were
obtained by Glasner, Megrelishvili and Uspenskij in
\cite{GM} and \cite{GMU}:

\begin{thm}
A compact metric dynamical system $(G,X)$ is tame if and only if
every element of $E(X)$ is a Baire 1 function from $X$ to itself.
\end{thm}

\begin{thm}
\label{th:main}
Let $X$ be a compact metric $G$-space.
The following conditions are equivalent:
\begin{enumerate}
\item
the dynamical system $(G,X)$ is hereditarily almost equicontinuous
(HAE);
\item
the dynamical system $(G,X)$ is
RN, that is,
admits a proper representation on a Radon--Nikod\'ym
Banach space;
\item
the enveloping semigroup $E(X)$ is metrizable.
\end{enumerate}
\end{thm}

For the definitions of HAE (hereditarily almost equicontinuous)
systems and the other undefined notions which appear in these theorems,
as well as for some further motivation and examples we refer the reader
to the papers \cite{GM}, \cite{G} and \cite{GMU}.

In \cite{G} I have shown that a {\em minimal} metrizable tame dynamical system
with a commutative acting group is PI and has zero topological entropy.
Recently Huang \cite{H}, and independently Kerr and Li \cite{KL}, improved
these results to show that under the same conditions a minimal tame system
is an almost 1-1 extension of its maximal equicontinuous factor
and is uniquely ergodic (see also \cite{HSY}).
In these works the authors make a heavy use of the structure
theory of minimal dynamical systems, as developed by
R. Ellis, W. Veech, Ellis-Glasner-Shapiro, McMahon and
van der Woude (see e.g. the survey \cite{G-lu}
and the references thereof).
However the main tool in both works (of Huang and Kerr and Li) is the
combinatorial notion of independence and the various related
notions of independence $n$-tuples.
In fact, Kerr and Li in their work \cite{KL}, use independence
to unify the theory of these various notions and in particular
they are able to characterize tame systems (which they call regular)
as those systems that (in some precise sense) do not admit infinite
independence sets (\cite[Proposition 6.4.2]{KL}).
In turn they use this characterization to define a notion of
relative regularity and develop the whole theory in the relative setup.

In the present work, which can be regarded as a continuation of
my work \cite{G}, I pursue purely structure theoretical
methods to recover the results of Huang and Kerr and Li
mentioned above, avoiding the combinatorial treatment altogether.
The key tool used in the proof, here as well as in \cite{G}, is a proposition
about diffused measures (Proposition \ref{diff} below),
which first appeared in \cite{G1}.

Section \ref{sur-sec} is a brief review of the structure theory
of minimal dynamical systems.
In Section \ref{so-sec}, I prove an analogue of an old theorem of
Ditor and Eifler \cite{DE}, which may have some independent interest.
It shows that when a continuous surjection $\pi: X \to Y$,
with $X$ and $Y$ compact metric, is semiopen then so is the induced
map $\pi_*: \Mf(X) \to \Mf(Y)$ on the spaces of probability measures
equipped with the weak$^*$ topology.
In Section \ref{diff-sec}, I pursue the idea of diffused measures,
first used in \cite{G1}, and prove the key Proposition \ref{diff}.
Section \ref{tame-sec} develops the theory of tame systems using
and extending results from \cite{G}.
In the final Section \ref{main-sec}, the main theorem is proved.

Except for the introductory Section \ref{sur-sec}, the group $\Ga$
is assumed to be Abelian.
For simplicity I handle only the case where the dynamical system is
metrizable and treat only the absolute (and not the more general
relative) case.

\section{A brief survey of abstract topological dynamics}\label{sur-sec}
This section is a brief review of the structure theory of
minimal dynamical systems. We will emphasize some
aspects which will be relevant in the present work.
For full details the reader is referred to
the books \cite{E2}, \cite{Gl1}, \cite{Au} and \cite{Vr}
and the review articles \cite{V} and \cite{G-lu}.

A {\em topological dynamical system} or briefly a system
is a pair $(X,\Ga)$, where
$X$ is a compact Hausdorff space and $\Ga$ an abstract infinite group
which acts on $X$ as a group of homeomorphisms.
A {\em sub-system} of \XGa is a closed invariant subset
$Y\subset X$ with the restricted action.
For a point $x\in X$, we let $\OGa (x)=\{\ga x:\ga \in \Ga\}$,
and $\OCGa (x)={\cls} \{\ga x:\ga \in \Ga\}$. These subsets of $X$
are called the {\em orbit} and
{\em orbit closure} of $x$ respectively.
We say that $(X,\Ga)$ is {\em point transitive} if
there exists a point $x\in X$ with a dense orbit.
In that case $x$ is called a {\em transitive point}.
If every point is transitive we say that $(X,\Ga)$ is a
{\em minimal system}. We say that $x\in X$ is an {\em
almost periodic} or a {\em minimal} point if $\OCGa(x)$
is a minimal system.

The dynamical system \XGa is {\em topologically transitive}
if for any two nonempty open subsets $U$ and $V$ of $X$
there exists some $\ga \in \Ga$ with $\ga U\cap V\ne\emptyset$.
Clearly a point transitive system is topologically transitive
and when $X$ is metrizable the converse holds as well:
in a metrizable topologically transitive system the
set of transitive points is a dense $G_\delta$ subset
of $X$.

The system \XGa is {\em weakly mixing} if the product
system $(X\times X, \Ga)$ (where $\ga(x,x')=(\ga x,\ga x'),
\ x,x'\in X,\ \ga \in \Ga$) is topologically transitive.

If $(Y,\Ga)$ is another system then a continuous onto
map $\pi:X\to Y$ satisfying $t\circ \pi=\pi\circ t$
for every $\ga \in \Ga$ is called a {\em homomorphism}
of dynamical systems. In this case we say that $(Y,\Ga)$
is a {\em factor} of $(X,\Ga)$ and also that
$(X,\Ga)$ is an {\em extension} of $(Y,\Ga)$.
With the system \XGa we associate the induced action
(the {\em hyper system} associated with \XGa)
on the compact space $2^X$ of closed subsets of $X$
equipped with the Vietoris topology.
A subsystem  $Y$ of $(2^X,\Ga)$ is a {\em
quasifactor} of \XGa if $\bigcup \{A: A\in Y\}=X$.

The system \XGa can always be considered as a quasifactor
of \XGa by identifying $x$ with $\{x\}$.
Recall that if \XtoYGa is a homomorphism then in general
$\pi^{-1}:Y\to 2^X$ is an upper-semi-continuous map
and that
$\pi:X\to Y$ is open iff $\pi^{-1}:Y\to 2^X$ is continuous,
iff $\{\pi^{-1}(y) : y\in Y\}$ is a quasifactor of \XGa.
When there is no room for confusion we write $X$ for
the system $(X,\Ga)$.

\br

We assume for simplicity that our acting group $\Ga$
is a discrete group. $\beta \Ga$ will denote the
Stone-\v {C}ech compactification of $\Ga$.
The universal properties of $\beta \Ga$ make it
\begin{itemize}
\item
a compact semigroup with right continuous
multiplication (for a fixed $p\in \beta \Ga$
the map $q\mapsto qp,\ q\in \beta \Ga$ is continuous),
and left continuous multiplication by
elements of $\Ga$, considered as elements of $\beta \Ga$
(for a fixed $\ga\in \Ga$
the map $q\mapsto \ga q,\ q\in \beta \Ga$ is continuous).
\item
a dynamical system $(\beta \Ga, \Ga)$ under
left multiplication by elements of $\Ga$.
\end{itemize}

The system $(\beta \Ga, \Ga)$ is the universal point
transitive $T$-system; i.e. for every point
transitive system \XGa and a point $x\in X$ with
dense orbit, there exists a homomorphism of
systems $(\beta \Ga, \Ga)\to (X,\Ga)$ which sends $e$,
the identity element of $\Ga$, onto $x$. For $p\in \beta \Ga$
we let $px$ denote the image of $p$ under this homomorphism.
This defines an ``action" of the semigroup $\beta \Ga$ on
every dynamical system. In fact, by universality there exists
a unique homomorphism $(\beta \Ga, \Ga)\to (E(X,\Ga),\Ga)$
onto the enveloping semigroup $E(X,\Ga)$ which is also a semigroup
homomorphism and we can interpret, and often do, the $\beta \Ga$
action on $X$ via this homomorphism.

When dealing with the hyper system $(2^X,\Ga)$ we write
$p\circ A$ for the image of the closed subset $A\subset X$
under $p\in \beta \Ga$ to distinguish it
from the (usually non-closed) subset $pA=\{px:x\in A\}$.
If $p$ is the limit of a net $\ga_i$ in
$\Ga$ then
$$
p \circ A =\{x\in X:
\text{there are a subnet $\ga_{i_j}$ and a net $x_j \in A$
with $x=\lim_j \ga_{i_j} x_j$}\}.
$$
We always have $pA\subset p\circ A$.

\br

The compact semigroup $\beta \Ga$ has a
rich algebraic structure.
For instance for countable $\Ga$
there are $2^c$ minimal
left (necessarily closed) ideals in
$\beta \Ga$ all isomorphic as systems
and each serving as a universal minimal system.
Each such minimal ideal, say $M$, has a subset $J$
of $2^c$ idempotents such that $\{vM:v\in J\}$ is a
partition of $M$ into disjoint isomorphic (non-closed)
subgroups. An idempotent in $\beta\Ga$ is called {\em minimal} if
it belongs to some minimal ideal.
A point $x$ in a dynamical system \XGa is
a minimal point iff there is some minimal idempotent $v$ in $\beta\Ga$
with $vx=x$, iff there exists some $v \in J$ with $vx=x$.

The group of dynamical system automorphisms of $(M,\Ga)$,\
$G={\Aut}(M,\Ga)$ can be identified with any one of the
groups $vM$ as follows:
with $\al\in vM$ we associate
the automorphism $\hat\al:(M,\Ga)\to (M,\Ga)$ given by
right multiplication $\hat\al(p)=p\al,\ p\in M$.
The group $G$ plays a central role in the algebraic theory.
It carries a natural $T_1$ compact topology, called by
Ellis the $\tau$-{\em topology}, which is weaker than the
relative topology induced on $G=uM$ as a subset of $M$.
The $\tau$-closure of a subset $A$ of $G$
consists of those $\beta\in G$ for which the set
${\gr}(\beta)=\{(p,p\beta):p\in M\}$
is a subset of the closure in $M\times M$
of the set $\bigcup\{{\gr}(\alpha):\alpha\in A\}$.
Both right and left multiplication on $G$ are $\tau$ continuous
and so is inversion.

It is convenient to fix a minimal left ideal $M$
in $\beta \Ga$ and an idempotent $u\in M$. As explained above we
identify $G$ with $uM$ and it follows that
for any subset $A\subset G$,
$$
{\cls}_\tau A = u(u\circ A) = G \cap (u\circ A).
$$
Also in this way we can consider
the ``action" of $G$ on every system \XGa via the
action of $\beta \Ga$ on $X$.
With every minimal
system \XGa and a point $x_0\in uX=\{x\in X: ux=x\}$ we
associate a $\tau$-closed subgroup
$$
\mathfrak{G}(X,x_0)=\{\alpha\in G:\alpha x_0=x_0\},
$$
the {\em Ellis group} of the pointed
system $(X,x_0)$. The quotient space $G/\mathfrak{G}(X,x_0)$
can be identified with the subset $uX \subset X$ via the map
$\al \mapsto \al x_0$
and the induced quotient $\tau$-topology is called the {\em
$\tau$-topology} on $uX$. Again the $\tau$-topology is weaker than
the relative topology induced on $uX$ as a subset of $X$, it is
$T_1$ and compact, and the closure operation is given by
$$
{\cls}_\tau A = u(u\circ A) = uX \cap (u\circ A),
\qquad A \subset uX.
$$

For a homomorphism $\pi: X \to Y$
with $\pi(x_0)=y_0$ we have
$$
\mathfrak{G}(X,x_0)\subset \mathfrak{G}(Y,y_0).
$$

For a $\tau$-closed subgroup $F$ of $G$ the {\em derived group}
$F'$ is given by:
\begin{equation}\label{prime}
F':=\bigcap\{{\cls}_\tau O : O\
\text{a $\tau$-open neighborhood of\ } u\  \text{in} \ F\}.
\end{equation}
$F'$ is a $\tau$-closed normal
(in fact characteristic) subgroup of $F$ and it is
characterized as the smallest $\tau$-closed subgroup
$H$ of $F$ such that $F/H$ is a compact
Hausdorff topological group.
In particular, for an Abelian $\Ga$, the topological group
$G/G'$ is the Bohr compactification of $\Ga$.

\br

A pair of points $(x,x')\in X\times X$ for a system \XGa
is called {\em proximal} if there exists a net
$\ga_i\in \Ga$ and a point $z\in X$ such that
$\lim t_ix=\lim t_ix'=z$
(iff there exists $p\in\beta \Ga$ with $px=px'$).
We denote by $P$ the set of proximal pairs in
$X\times X$.
We have
$$
P=\bigcap\ \{\Ga V:V\
{\text{\rm a neighborhood of the diagonal in}}\ X\times X\}.
$$
A system \XGa is called {\em proximal} when $P=X\times X$
and {\em distal} when $P=\Delta$, the diagonal in $X\times X$.
It is called {\em strongly proximal} when the following
much stronger condition holds: the dynamical system
$(\Mf(X),\Ga)$, induced on the compact space $\Mf(X)$ of probability
measures on $X$, is proximal.
A minimal system \XGa is called {\em point distal}
if there exists a point $x\in X$ such that if
$x,x'$ is a proximal pair then $x=x'$.

The {\em regionally proximal relation\/} on
$X$ is defined by
$$
Q=\bigcap\ \{\ol {\Ga V}:V\
{\text{\rm a neighborhood of $\Delta$ in}}\ X\times X\}.
$$
It is easy to verify that $Q$ is trivial --- i.e. equals
$\Delta$ --- iff the system is equicontinuous.

\br

An extension \XtoYGa of minimal systems
is called a {\em proximal extension}
if the relation $R_\pi=\{(x,x'):\pi(x)=\pi(x')\}$
satisfies $R_\pi\subset P$ and
a {\em distal extension} when $R_\pi\cap P=\Del$.
One can show that every distal extension is open.
$\pi$ is a {\em highly proximal (HP) extension} if for every closed
subset $A$ of $X$ with $\pi(A)=Y$, necessarily $A=X$.
It is easy to see that a HP extension is proximal.
In the metric case an extension \XtoYGa of minimal systems is
HP iff it is an {\em almost  1-1 extension}, that is the set
$\{y\in Y:  {\text{ with $\pi^{-1}(y)$ is a singleton}}\}$
is a dense $G_\del$ subset of $Y$.
The map $\pi$ is {\em strongly proximal} if for every $y\in Y$
and every probability measure $\nu$ with ${\supp}\nu \subset \pi^{-1}(y)$,
there exists a net $\ga_i\in \Ga$ and a point $x\in X$ such that
$\lim_i \ga_i \nu =\delta_x$ in the weak$^*$ topology on the space
$\Mf(X)$ of probability measures on $X$.
The extension $\pi$ is called an
{\em equicontinuous extension} if for every $\ep$,
a neighborhood of the diagonal
$\Delta = \{(x,x): x\in X\} \subset X\times X$,
there exists a neighborhood of the diagonal $\del$ such that
$\ga (\del\cap R_\pi) \subset \ep$ for every $\ga \in \Ga$.
In the metric case an equicontinuous extension is also called
an {\em isometric extension}.
The extension $\pi$ is a {\em weakly mixing extension}
when $R_\pi$ as a subsystem of the product system
$(X\times X,\Ga)$ is topologically transitive.

\br

The algebraic language is particularly suitable for dealing
with such notions. For example an extension
\XtoYGa of minimal systems is a proximal extension iff
the Ellis groups
$\mathfrak{G} (X,x_0)=A$ and $\mathfrak{G} (Y,y_0)=F$ coincide.
It is distal iff for every $y\in Y$,
and $x\in \pi^{-1}(y),\
\pi^{-1}(y)=\mathfrak{G}(Y,y)x$;\
iff:
\begin{quote}
for every $y=py_0\in Y$, $p$ an element of $M$,
$\pi^{-1}(y)=p\pi^{-1}(y_0)=
pFx_0$, where $F=\mathfrak{G}(Y,y_0)$.
\end{quote}
In particular $(X,\Ga)$ is distal iff
$Gx=X$ for some (hence every) $x\in X$.
The extension $\pi$ is an equicontinuous extension
iff it is a distal extension and, denoting
$\mathfrak{G} (X,x_0)=A$ and $\mathfrak{G} (Y,y_0)=F$,
$$
F'\subset A.
$$
In this case, setting $A_0=\bigcap_{g\in F} gAg^{-1}$,
the group $F/A_0$ is the
group of the {\em group extension\/}\ $\tilde \pi$ associated with
the equicontinuous extension $\pi$.
More precisely, there exists a minimal
dynamical system $(\tilde X,\Ga)$, with
${\mathfrak{G}}(\tilde{X},\tilde{x}_0)=A_0$,
on which the compact Hausdorff topological
group $K=F/A_0$ acts as a group of automorphisms and
we have the following commutative diagram
\begin{equation*}
\xymatrix
{
\tilde{X} \ar[dd]_{\tilde{\pi}}\ar[dr]^{\phi}  & \\
 & X\ar[dl]^{\pi}\\
Y &
}
\end{equation*}
where $\tilde\pi: \tilde X \to Y\cong \tilde X/ K$ is
a group extension and so is the extension
$\phi: \tilde X \to X\cong \tilde X/L$
with $L=A/A_0 \subset F/A_0 = K$. ($\tilde X =X$ iff
$A$ is a normal subgroup of $F$.)

\br

A minimal system \XGa is called {\em incontractible}
if the union of minimal subsets is dense in every
product system $(X^n, \Ga)$. This is the case iff
$p\circ Gx=X$ for some (hence every) $x\in X$
and $p\in M$. When $\Ga$ is Abelian $Gx$ is always dense
in $X$ so that every minimal system is incontractible.
However the following relative notion is an important tool
even when $\Ga$ is Abelian.

We say that \XtoYGa is a RIC ({\em relatively incontractible})
{\em extension} if:
\begin{quote}
for every $y=py_0\in Y$, $p$ an element of $M$,
$\pi^{-1}(y)=p\circ u\pi^{-1}(y_0)=
p\circ Fx_0$, where $F=\mathfrak{G}(Y,y_0)$.
\end{quote}
One can show that every RIC extension is open and that
every distal extension is RIC. It then
follows that every distal extension is open.

We have the following theorem from \cite{EGS}
about the interpolation of equicontinuous extensions.
For a proof see \cite{G1}, Theorem X.2.1.

\begin{thm}\label{eq-ext}
Let $\pi: X \to Y$ be a RIC extension of minimal systems.
Fix a point $x_0\in X$ with $ux_0=x_0$ and let $y_0=\pi(x_0)$.
Let $A=\mathfrak{G}(X,x_0)$ and $F=\mathfrak{G}(Y, y_0)$.
Then there exists a commutative diagram of pointed systems
\begin{equation*}
\xymatrix
{
(X,x_0) \ar[dd]_{\pi}\ar[dr]^{\sig}  & \\
 & (Z,z_0)\ar[dl]^{\rho}\\
(Y,y_0) &
}
\end{equation*}
such that  $\rho$ is an equicontinuous extension with
Ellis group $\mathfrak{G}(Z,z_0)=AF'$ and
the extension $\rho$ is an isomorphism iff $AF'=F$.
Moreover if
\begin{equation*}
\xymatrix
{
(X,x_0) \ar[dd]_{\pi}\ar[dr]^{\sig'}  & \\
 & (Z',z'_0)\ar[dl]^{\rho'}\\
(Y,y_0) &
}
\end{equation*}
is another such diagram with $\rho'$ an equicontinuous
extension then there exists a homomorphism $(Z,z_0) \to (Z',z'_0)$.
\end{thm}

\br

Given a homomorphism $\pi:(X,\Ga) \to (Y,\Ga)$ of minimal metric systems,
there are several standard constructions of associated ``shadow diagrams".
In the {\em O shadow diagram}
\begin{equation*}
\xymatrix
{
X \ar[d]_{\pi}  & X^* =X \vee Y^*\ar[l]_-{\tet^*} \ar[d]^{\pi^*} \\
Y  & Y^*\ar[l]^{\tet}
}
\end{equation*}
the map $\pi^*$ is open
and the maps $\theta$ and $\theta^*$ are almost 1-1.
The explicit constructions is as follows.
The set valued map $\pi^{-1}: Y \to 2^X$ (where the latter is
the compact space of closed subsets of $X$, equipped with the
Hausdorff, or Vietoris, topology) is uppersemicontinuous
and we let $Y_0\subset Y$ be the set of continuity points
of this map. Set
$Y^*={\cls}\{\pi^{-1}(y): y\in Y_0\}\subset 2^X$,
and $X^*=X\vee Y^* = {\cls}\{(x,\pi^{-1}(y)): y\in Y_0,\ \pi(x)=y\}
\subset X\times Y^*$.
By the uppersemicontinuity of $\pi^{-1}$ every $y^*\in Y^*$
is contained in a fiber $\pi^{-1}(y)$ for some $y\in Y$
and we let $\tet(y^*)=y$.
The maps $\pi^*$ and $\tet^*$ are the restriction to $X^*$
of the coordinate projections on $X$ and $Y^*$ respectively.
One then shows that
$X^*=\{(x,y^*): x \in y^*\in Y^*\}$
and that indeed, $\pi^*$ is open
and the maps $\theta$ and $\theta^*$ are highly proximal.
The O shadow diagram collapses, i.e. $Y=Y^*$, $X=X^*$
and $\pi=\pi^*$ iff $\pi: X\to Y$ is an open map; iff
the map $\pi^{-1}: Y \to 2^X$ is continuous.

In the {\em RIC-shadow diagram}
\begin{equation*}
\xymatrix
{
X \ar[d]_{\pi}  & X^*=X\vee Y^* \ar[l]_-{\tet^*} \ar[d]^{\pi^*} \\
Y  & Y^*\ar[l]^{\tet}
}
\end{equation*}
$ \pi^*$ is RIC and $\theta, \theta^*$
are proximal (thus we still have
$A={\mathfrak{G}}(X,x_0)={\mathfrak{G}}( X^*, x^*_0)$ and
$F={\mathfrak{G}}(Y,y_0)= {\mathfrak{G}}(Y^*, y^*_0)$).
The concrete description of these objects uses
quasifactors and the circle operation:
$$
Y^*= \{p\circ Fx_0:p\in M\}\subset  2^X,\qquad
X^*= \{(x, y^*):x\in  y^*\in Y^*\}
\subset X\times  Y^*
$$
and
$$
\theta(p\circ Fx_0)=py_0, \quad
\theta^*(x, y^*)=x,\quad
\pi^*(x, y^*)= y^*,\quad (p\in M),
$$
where $F={\mathfrak{G}}(Y,y_0)$.
The map $\theta$ is an isomorphisms (hence $\pi=\pi^*$)
when and only when $\pi$ is already RIC.

Finally we say that $\pi:(X,\Ga) \to (Y,\Ga)$ has a
{\em relatively invariant measure} (RIM),
if there exists a projection $P:C(X) \to C(Y)$ such that
\begin{enumerate}
\item
$P(f) \ge 0$ for $f\ge 0$ in $C(X)$.
\item
$P(\ch)=1.$
\item
$P(h \circ \pi)= h$ for every $h \in C(Y)$.
\item
$P(f\circ \ga) = P(f) \circ \ga$ for every $f\in C(X)$ and $\ga\in \Ga$.
\end{enumerate}
This property is equivalent to the existence of a continuous {\em section},
i.e. a continuous $\Ga$ equivariant map $y \mapsto \la_y$
from $Y$ into $\Mf(X)$ such that $\pi(\la_y)=\del_y$ for every $y\in Y$.
Here and in the sequel we use the same letter $\pi$ to denote the induced
map $\pi: \Mf(X) \to \Mf(Y)$ on the spaces of probability measures.
Sometimes though we will write $\pi_*$ for the induced map.

In the {\em RIM shadow diagram}
\begin{equation*}
\xymatrix
{
X \ar[d]_{\pi}  & \tilde{X}=X \vee \tilde{Y}
 \ar[l]_-{\tilde{\tet}} \ar[d]^{\tilde{\pi}} \\
Y  & \tilde{Y}\ar[l]^{\tet}
}
\end{equation*}
the map $\tilde \pi$ has a RIM
and the maps $\theta$ and $\tilde\theta$ are strongly proximal.
It can be shown that every isometric extension has a RIM and is open.
See \cite{G0} for more details, also a treatment of SPI systems can
be found in \cite{G1}.

\br

We say that a minimal system \XGa is a
{\em strictly PI system} if there is an ordinal $\eta$
(which is countable when $X$ is metrizable)
and a family of systems
$\{(W_\iota,w_\iota)\}_{\iota\le\eta}$
such that (i) $W_0$ is the trivial system,
(ii) for every $\iota<\eta$ there exists a homomorphism
$\phi_\iota:W_{\iota+1}\to W_\iota$ which is
either proximal or equicontinuous
(isometric when $X$ is metrizable), (iii) for a
limit ordinal $\nu\le\eta$ the system $W_\nu$
is the inverse limit of the systems
$\{W_\iota\}_{\iota<\nu}$,  and
(iv) $W_\eta=X$.
We say that \XGa is a {\em PI-system} if there
exists a strictly PI system $\tilde X$ and a
proximal homomorphism $\theta:\tilde X\to X$.

If in the definition of PI-systems we replace
proximal extensions by HP extensions (almost 1-1 extensions
in the metric case) we get the notion of {\em HPI} ({\em AI-systems}
in the metric case).
If we replace the proximal extensions by trivial
extensions (i.e.\ we do not allow proximal
extensions at all) we have {\em I-systems}.
In this terminology the structure theorem for distal
systems (Furstenberg \cite{F}, 1963) can be stated as follows:
\begin{thm}\label{FST}
A metric minimal system is distal iff it is an I-system.
\end{thm}

And the Veech-Ellis structure theorem for point distal systems
(Veech \cite{V}, 1970 and Ellis \cite{E1}, 1973).

\begin{thm}
A metric minimal dynamical system is point distal iff it is an
AI-system.
\end{thm}

The structure theorem for the general minimal system is proved in
\cite{EGS} and \cite{McM} (see also \cite{V}) and asserts that
every minimal system admits a canonically defined proximal
extension which is a weakly mixing RIC extension of a strictly PI system.
Both the Furstenberg and the Veech-Ellis structure theorems are
corollaries of this general structure theorem.

\begin{thm}[Structure theorem for minimal systems]\label{structure}
Given a minimal system $(X,\Ga)$, there exists an ordinal $\eta$
(countable when $X$ is metrizable) and a canonically defined
commutative diagram (the canonical PI-Tower)
\begin{equation*}
\xymatrix
        {X \ar[d]_{\pi}             &
     X_0 \ar[l]_{{\theta}^*_0}
         \ar[d]_{\pi_0}
         \ar[dr]^{\sigma_1}         & &
     X_1 \ar[ll]_{{\theta}^*_1}
         \ar[d]_{\pi_1}
         \ar@{}[r]|{\cdots}         &
     X_{\nu}
         \ar[d]_{\pi_{\nu}}
         \ar[dr]^{\sigma_{\nu+1}}       & &
     X_{\nu+1}
         \ar[d]_{\pi_{\nu+1}}
         \ar[ll]_{{\theta}^*_{\nu+1}}
         \ar@{}[r]|{\cdots}         &
     X_{\eta}=X_{\infty}
         \ar[d]_{\pi_{\infty}}          \\
        pt                  &
     Y_0 \ar[l]^{\theta_0}          &
     Z_1 \ar[l]^{\rho_1}            &
     Y_1 \ar[l]^{\theta_1}
         \ar@{}[r]|{\cdots}         &
     Y_{\nu}                &
     Z_{\nu+1}
         \ar[l]^{\rho_{\nu+1}}          &
     Y_{\nu+1}
         \ar[l]^{\theta_{\nu+1}}
         \ar@{}[r]|{\cdots}         &
     Y_{\eta}=Y_{\infty}
    }
\end{equation*}
where for each $\nu\le\eta, \pi_{\nu}$
is RIC, $\rho_{\nu}$ is isometric, $\theta_{\nu},
{\theta}^*_{\nu}$ are proximal and
$\pi_{\infty}$ is RIC and weakly mixing.
For a limit ordinal
$\nu ,\  X_{\nu}, Y_{\nu}, \pi_{\nu}$
etc. are the inverse limits (or joins) of
$ X_{\iota}, Y_{\iota}, \pi_{\iota}$ etc. for $\iota
< \nu$.
Thus $X_\infty$ is a proximal extension of $X$ and a RIC
weakly mixing extension of the strictly PI-system $Y_\infty$.
The homomorphism $\pi_\infty$ is an isomorphism (so that
$X_\infty=Y_\infty$) iff $X$ is a PI-system.
\end{thm}

Two further corollaries of this theorem are
the theorems of Bronstein on the structure of PI systems, \cite{Bro}
and of van der Woude on HPI systems, \cite{vdW}.
Here we will use the latter which I now proceed to describe.

A homomorphism $\pi:(X,\Ga) \to (Y,\Ga)$ is called {\em semiopen}
if the interior of $\pi(U)$ is nonempty for every nonempty open subset
$U$ of $X$. When $X$ is minimal every
$\pi:(X,\Ga) \to (Y,\Ga)$ is semiopen.
We will say that a subset $W \subset X\times X$ is a {\em S-set}
if it is closed invariant topologically transitive and the restriction
to $W$ of the projection maps are semiopen.

\begin{thm}[van der Woude]\label{vdw}
A minimal system $(X,\Ga)$ is HPI iff every S-set in $X\times X$
is minimal.
\end{thm}

\br

\section{On semiopen maps}\label{so-sec}

A result of Ditor and Eifler from 1972, \cite{DE} asserts that a
continuous surjection $\pi: X\to Y$ between compact Hausdorff
spaces $X$ and $Y$ is open iff the induced map $\pi_*:
\Mf(X) \to \Mf(Y)$ is an open surjection.
In the course of the proof of our main theorem (Theorem \ref{almost-auto})
we will need an analogous result (in the metric case) for semiopen maps
(Theorem \ref{so-thm}), which we now proceed to establish.
First we need two preliminary lemmas.

\begin{lem}\label{so-lem}
Let $\pi: X\to Y$ be a continuous surjection between
compact Hausdorff spaces. The conditions {\rm 1} and {\rm 2}
below are equivalent.
If $X$ is metrizable then the three conditions
are equivalent:
\begin{enumerate}
\item
$\pi$ is semiopen.
\item
The preimage of every dense subset in $Y$ is dense in $X$.
\item
The set
\begin{equation*}
X_0 = \{x\in X: \text{ the set valued map $\pi^{-1}: Y \to 2^X$
is continuous at $\pi(x)$}\}
%=\{x\in X: \text{ the  map $\pi: X \to Y$ is open at $x$}\}
\end{equation*}
is dense in $X$.
\end{enumerate}
\end{lem}

\begin{proof}
The equivalence of 1 and 2 is straightforward. For any
continuous surjection $\pi: X\to Y$ the corresponding set map
$\pi^{-1}: Y \to 2^X$ is uppersemicontinuous and, when $X$ is metrizable,
this implies that it has a dense $G_\del$ subset $Y_0\subset Y$ of
continuity points.
Assuming 2 we conclude that $X_0=\pi^{-1}(Y_0)$ is a dense $G_\del$
subset of $X$. Conversely if 3 is valid and $U\subset X$ is
open and nonempty, then $U \cap X_0 \not=\emptyset$ and if $x_0$
is any point in this intersection then $\pi$ is open at $x_0$,
so that $\pi(U)$ is a neighborhood of $\pi(x_0)$ and we conclude
that $\inte (\pi(U))\not=\emptyset$.
\end{proof}

\begin{lem}\label{f}
Let $\pi: X\to Y$ be a continuous surjection between
compact metric spaces. Let $f: X \to \R$
be a continuous function and define $f^*: Y \to \R$
by
$$
f^*(y)=\sup\{f(x): \pi(x)=y\}.
$$
Then $f^*$ is continuous at every point of the set
$$
Y_0=\{y\in Y: \text{ the set valued map $\pi^{-1}: Y \to 2^X$
is continuous at $y$}\}.
$$
\end{lem}

\begin{proof}
Fix $y\in Y_0$ and suppose $y_n \to y$ is a convergent sequence.
For each $n$ let $x_n \in X$ satisfy $\pi(x_n)=y_n$ and $f^*(y_n)=f(x_n)$.
Let $x_{n_k} \to x$ be a convergent subsequence. Then $\pi(x)=y$ hence
$$
f^*(y) \ge f(x)=\lim_{k\to\infty}f(x_{n_k}) = \lim_{k\to\infty}f^*(y_{n_k}).
$$
Since this is true for every partial limit of $f^*(y_n)$ we conclude that
$$
f^*(y) \ge \limsup_{n\to\infty}f^*(y_n).
$$
On the other hand if $f^*(y)=f(x)$ with $\pi(x)=y$ then, since $\pi$
is open at $x$, we can find a sequence $x_n$ with $\pi(x_n)=y_n$
and $x_n \to x$, so that
$$
f^*(y)=f(x) = \lim_{n\to\infty}f(x_n) \le
\liminf_{n\to\infty}f^*(y_n).
$$
\end{proof}

\begin{thm}\label{so-thm}
Let $\pi: X \to Y$ be a continuous surjection between compact
metric spaces which is semiopen. Then the induced map $\pi_*:
\Mf(X) \to \Mf(Y)$ is a semiopen surjection.
\end{thm}

\begin{proof}
Let
$$
Y_0=\{y\in Y: \text{ the set valued map $\pi^{-1}: Y \to 2^X$
is continuous at $y$}\}
$$
and
$$
X_0 =\{x\in X: \text{ the set valued map $\pi^{-1}: Y \to 2^X$
is continuous at $\pi(x)$}\}.
$$
Then $X_0=\pi^{-1}(Y_0)$ and by Lemma \ref{so-lem}
$Y_0$ and $X_0$ are dense $G_\del$ subsets of $Y$ and $X$
respectively.
Let $\Mf_0(X)$ be the collection of measures of the form
$\sum_{i=1}^n c_i \del_{x_i}$, where $0 \le c_i \le 1$,
$\sum_{i=1}^n c_i=1$.
If in addition we require that each $x_i$ is in $X_0$
we obtain the smaller collection $\Mf_{00}(X)$. Clearly
$\Mf_{00}(X)$, and hence also $\Mf_0(X)$, are dense in $\Mf(X)$.
We define $\Mf_0(Y)$ analogously as the collection of measures in $\Mf(Y)$
of the form $\nu= \sum_{i=1}^n c_i \del_{y_i}$.
Again $\Mf_0(Y)$ is dense in $\Mf(Y)$.

Let $\Ucal \subset \Mf(X)$ be a closed set with nonempty interior.
We have to show that $\inte (\pi_*(\Ucal))\not=\emptyset$.
Suppose to the contrary that $\inte (\pi_*(\Ucal))=\emptyset$.
Fix $\mu_0 =\sum_{i=1}^m c_i \del_{x_i} \in \Mf_{00}(X) \cap \inte(\Ucal)$
and let $\nu_0=\pi_*(\mu_0) = \sum_{i=1}^m c_i \del_{y_i}$,
where $y_i=\pi(x_i)$.
Let $\nu_j \in \Mf_0(Y)\setminus \pi_*(\Ucal)$ be a sequence which converges
to $\nu_0$.
Set $Q_j =\pi_*^{-1}(\nu_j)$. Each $Q_j$ is a closed and convex subset
of $\Mf(X)$ and with no loss of generality we assume that
$Q=\lim_{j\to\infty} Q_j$ exists in $2^{\Mf(X)}$. Then $Q$ is
a compact convex subset of $\Mf(X)$ with $\pi_*(\tet)=\nu_0$
for every $\tet \in Q$.

If $\mu_0 \in Q$ then eventually $Q_j \cap \Ucal \not=\emptyset$,
hence $\nu_j \in \pi_*(\Ucal)$, contradicting our assumption. Thus
we have $\mu_0 \not\in Q$ and by the separation theorem there
exist a function $f\in C(X)$ and $\ep>0$ with
\begin{equation}\label{sep}
\mu_0(f) \ge \tet(f) +\ep  \  {\text{for every $\tet\in Q$}}.
\end{equation}
Define the associated function
$$
f^*(y)=\sup\{f(x): \pi(x)=y\}.
$$
Each measure $\nu_j$ has the form
$\sum_{i=1}^{m_j} c_{j,i} \del_{y_{j,i}}$ and we choose points
$x_{j,i}$ with $\pi(x_{j,i})=y_{j,i}$ such that
$f^*(y_{j,i}) = f(x_{j,i})$. Now form the measures
$\mu_j = \sum_{i=1}^{m_j} c_{j,i} \del_{x_{j,i}}$ and
assume, with no loss of generality, that
$\mu =\lim_{j\to\infty}\mu_j$ exists in $\Mf(X)$.
Since $\mu_j\in Q_j$ for each $j$, we have $\mu \in Q$.

Note that by our construction $\nu_j(f^*) = \mu_j (f)$ for every $j$.
By assumption the set ${\supp}\nu_0 =\{y_1,\dots,y_m\}$
is a subset of $Y_0$ and therefore, by Lemma \ref{f},
each $y_i$ is a continuity point for $f^*$.
From this fact it is easy to deduce that
$$
\lim_{j\to \infty}\nu_j(f^*) = \nu_0(f^*).
$$
It then follows that
$$
\mu(f) = \lim_{j\to \infty}\mu_j(f) =
\lim_{j\to \infty}\nu_j(f^*) = \nu_0(f^*)=
\sum_{i=1}^{m} c_i f^*(y_i) = \mu_0 (f^* \circ \pi)
\ge \mu_0(f).
$$
This contradicts (\ref{sep}) and our proof is complete.
\end{proof}

Recall the following well known result;
for completeness we include a proof.

\begin{lem}\label{so-min}
Let $\pi: (X,\Ga) \to (Y,\Ga)$ be a homomorphism between minimal
systems. Then  $\pi$ is semiopen.
\end{lem}

\begin{proof}
Let $W\subset X$ be a closed set with nonempty interior.
By minimality of \XGa there is a finite set $\{\ga_1,\dots,\ga_n\}
\subset \Ga$ with $X = \bigcup_{i=1}^n \ga_i W$.
Therefore $Y = \bigcup_{i=1}^n \pi(\ga_i W)$
and it follows that for some $i$ the interior of the closed
set $\pi(\ga_i W) = \ga_i\pi(W)$ is nonempty. Thus, as required, also
$\inte (\pi(W))\not=\emptyset$.
\end{proof}

\br

\section{A key proposition on diffused measures}\label{diff-sec}
As explained above we fix a minimal ideal $M$ in $\beta \Ga$ and let
$u$ be an idempotent in $M$.
We denote the subgroup $uM$ of $M$ by $G$ and identify it
with the group of automorphisms of the universal $\Ga$-minimal
system $(M,\Ga)$, where for $\al\in G$ the corresponding
automorphism $R_\al: M \to M$ is given by right multiplication
$p \mapsto p \al$.
For an Abelian $\Ga$ each subgroup $vM\subset M$, where $v$ is an idempotent
in $M$, is dense in $M$ and it follows that the $G$-dynamical
system $(M,G)$ (where $G$ acts by right multiplication) is minimal.
From now on we always assume that $\Ga$ is an Abelian group.

\begin{lem}\label{sat}
Let $(Y^*,\Ga)$ be a metric minimal system and
$\tet:(Y^*,\Ga) \to (Y,\Ga)$ its maximal equicontinuous factor.
Suppose further that $\tet$ is an almost 1-1 extension.
Let $O$ be a nonempty open subset of $Y^*$.
\begin{enumerate}
\item
There is a nonempty open subset $V\subset O$ such that
$\tet^{-1}(\tet(V))\subset O$.
\item
There is a nonempty open subset $W\subset O$ such that
${\cls}_\tau(W \cap uY^*) \subset O$.
\end{enumerate}
\end{lem}

\begin{proof}

1.\
Suppose $O \subset Y^*$ is a nonempty open set for which the statement of
the lemma fails. Choose a point $y^* \in O$ such that
$\tet^{-1}(\tet(y^*))=\{y^*\}$ and let $V_n\subset O$ be a sequence of
open balls centered at $y^*$ with ${\diam}(V_n) \searrow 0$.
By assumption there are pairs of points $z_n \in V_n$
and $z'_n\not\in O$ with $\tet(z_n)=\tet(z'_n)$. However,
as $\lim_{n\to\infty} \tet^{-1}(\tet(z_n)) = \{y^*\}$,
we have $z'_n \to y^*$ in contradiction of the fact that $O$
is a neighborhood of $y^*$.

2.\
Since $Y$ is equicontinuous its $\tau$-topology coincides with its
compact Hausdorff group topology.
Since $\tet$ is an almost 1-1 map,
the restriction $\tet\rest uY^* : uY^* \to Y$ is a
homeomorphism of $uY^*$, equipped with the $\tau$ topology, onto $Y$.
Let $O \subset Y^*$ be a nonempty open set. Let $V \subset O$
be as in part 1, and let $W$ be a
nonempty open subset such that $\ol{W} \subset V$.
Now
\begin{align*}
{\cls}_\tau (W \cap uY^*) & = \tet^{-1}(\ol{\tet(W)}) \cap uY^*\\
& \subset  \tet^{-1}(\ol{\tet(W)})
 = \tet^{-1}(\tet(\ol{W}))\\
& \subset  \tet^{-1}(\tet(V)) \subset O.
\end{align*}
\end{proof}

\begin{lem}\label{tau}
Let \XGa be a minimal metric system and let
$\phi:(X,\Ga) \to (Y,\Ga)$ be its maximal equicontinuous factor.
Suppose further that $Y$ is infinite and that we have the following diagram
$$
(X,\Ga) \overset{\pi}{\to} (Y^*,\Ga) \overset{\tet}{\to} (Y,\Ga),
$$
where $\pi$ is an isometric extension, $\tet$ is an almost
1-1 extension and $\phi= \tet \circ \pi$.
Let $U$ be a nonempty open subset of $X$.
Fix a minimal ideal $M \subset \beta\Ga$ and an idempotent $u\in M$.
Then
$$
{\cls}_\tau (U\cap uX) \supset \pi^{-1}(\pi(U)) \cap uX.
$$
\end{lem}

\begin{proof}
Fix a point $x_0 \in X$ with $ux_0=x_0$ and let $y_0$ and $y^*_0$
be its images in $Y$ and $Y^*$ respectively.
Let
\begin{gather*}
A=\mathfrak{G}(X,x_0)=\{\alpha\in G:\alpha x_0=x_0\}, \ \text{and}\ \\
F = \mathfrak{G}(Y,y_0)=
\mathfrak{G}(Y^*,y^*_0)=\{\alpha\in G:\alpha y_0=y_0\}=
\{\alpha\in G:\alpha y^*_0=y^*_0\}.
\end{gather*}
The assumption that $\pi$ is an isometric extension implies that
$B = F/A$ is a homogeneous space of a Hausdorff compact topological group.
The fact that \YGa is the maximal equicontinuous factor of \XGa
implies that $F \supset G'$ and that $G'A=AG'=F$.
Let $\sig: M \to X$ denote the evaluation map $p \mapsto px_0$.

Let $U$ be a nonempty open subset of $X$.
Set $\tilde{U}=\sig^{-1}(U)=
\{p\in M: px_0 \in U\}$.
Then $\tilde{U}$ is a nonempty open subset of $M$ and
by minimality of the $G$-system $(M,G)$,
the collection $\{\tilde{U}\al: \al \in G\}$ is an open cover of $M$.
Choose a finite subcover, say $\{\tilde{U}\al_i: i=1,2,\dots,n\}$.
Now
$$
\bigcup_{i=1}^n  {\cls}_\tau(\tilde{U}\al_i \cap G)
= \bigcup_{i=1}^n {\cls}_\tau(\tilde{U} \cap G)\al_i = G;
$$
hence ${\cls}_\tau(\tilde{U}\cap G)$ has a nonempty $\tau$-interior.
Since ${\cls}_\tau(\tilde{U} \cap G)$ is also $\tau$-closed, it must
contain a left translate of $G'$, say $\beta G'$ for some $\beta \in G$
(this follows from the definition of $G'$, see equation (\ref{prime})).
Projecting back to $X$ via $\sig$ we get
\begin{align*}
{\cls}_\tau (U \cap uX ) & =
\sig({\cls}_\tau (\tilde{U}\cap G)) \\
&\supset
\beta G' x_0 = \beta G' A x_0 =
\beta Fx_0 \\
& = (\beta F \beta^{-1})\beta x_0 = \pi^{-1}(\pi(\beta x_0)).
\end{align*}
As $(Y,\Ga)$ is equicontinuous and since $\tet$ is an almost 1-1 map,
it follows that the restriction $\tet\rest uY^* : uY^* \to Y$ is a
homeomorphism of $uY^*$ equipped with the $\tau$ topology
onto $Y$ equipped with its compact Hausdorff group topology.

Let $O=\pi(U)$, then $O$ is a nonempty open subset of $Y^*$
and by Lemma \ref{sat} there is a nonempty open subset $W\subset O$
such that ${\cls}_\tau(W \cap uY^*) \subset O$. Set $U_1=\pi^{-1}(W)
\cap U$. Then $U_1$ is a nonempty open subset of $X$ and by the above argument
there exists $\beta_1\in G$ with $\pi^{-1}(\pi(\beta_1 x_0))
\subset {\cls}_\tau (U_1 \cap uX )\subset {\cls}_\tau (U \cap uX )$. Now
\begin{align*}
\pi(\beta_1 x_0) &
\in \pi({\cls}_\tau (U_1 \cap uX ))\\
& = {\cls}_\tau (\pi(U_1) \cap uY^*)\\
& ={\cls}_\tau (W \cap uY^* ) \subset O = \pi(U).
\end{align*}

Thus we have shown that for every nonempty open subset $U\subset X$,
the set ${\cls}_\tau (U \cap uX )$ contains a full fiber
$\pi^{-1}(y^*)$ for some $y^*\in \pi(U)$. Since $\pi$ is an open map
we conclude that
$$
{\cls}_\tau (U\cap uX) \supset \pi^{-1}(\pi(U)) \cap uX
$$
as required.
\end{proof}

\begin{prop}\label{diff}
Let \XGa be a minimal metric system and let
$\phi:(X,\Ga) \to (Y,\Ga)$ be its maximal equicontinuous factor.
Suppose further that $Y$ is infinite and that we have the following diagram
$$
(X,\Ga) \overset{\pi}{\to} (Y^*,\Ga) \overset{\tet}{\to} (Y,\Ga),
$$
where $\pi$ is an isometric extension, $\tet$ is an almost
1-1 extension and $\phi= \tet \circ \pi$.
For every $y^*\in Y^*$ the fiber $\pi^{-1}(y^*)$ has the structure of a
homogeneous space of a compact Hausdorff topological group
and we let $\la_{y^*}$ be the corresponding Haar measure on
this fiber.
Thus $\pi$ is a RIM and open extension and
$y^*\mapsto \la_{y^*}$, \ $Y^* \to \Mf(X)$, is the corresponding section.
Let $\La: \Mf(Y^*) \to \Mf(X)$, defined by
\begin{equation*}
\La(\nu) =\int_{Y^*} \la_{y^*}\, d\nu(y^*),
\end{equation*}
be the associated affine injection.
Set $\Mf_m(X)=\{\La(\nu): \nu \in \Mf(Y^*)\}$.
Then the set
$$
R=\{\nu\in \Mf(X): \text{the orbit closure of $\nu$
meets $\Mf_m(X)$}\}
$$
is a dense $G_\del$ subset of $\Mf(X)$.
\end{prop}

\begin{proof}
Fix a compatible metric $d$ on $\Mf(X)$ and
let $\kappa \in \Mf(X)$ and $\ep, \eta >0$ be given.
Find an atomic measure $\la=\frac{1}{n}\sum_{i=1}^n \del_{x_i}$
such that $d(\ka,\la) < \ep/2$.
Since $Y$ is infinite and $\pi$ is an open map we can also assume that
$\pi(x_i)\not =\pi(x_j)$ for $i \not = j$.
Choose open disjoint neighborhoods $U_i$ of $x_i$,
so small that every measure of the form
$\mu=\frac{1}{n}\sum_{i=1}^n \mu_i,\ \mu_i \in \Mf(X)$
 with ${\supp}\mu_i \subset U_i$,
will satisfy $d(\mu,\la) < \ep/2$, and hence also
$d(\mu,\ka) < \ep$.

Set $\nu=\frac{1}{n}\sum_{i=1}^n \la_{y^*_i}$ with $y^*_i=\pi(x_i)$.
For each $y^*_i$ choose points $\{x'_{i,j}\}_{j=1}^k
\in \pi^{-1}(y^*_i)$ so that $d(\mu',\nu) < \ep/2$,
where
$$
\mu'=\frac{1}{nk}\sum_{i=1}^n \sum_{j=1}^k \del_{x'_{i,j}}.
$$

By Lemma \ref{tau}
\begin{align*}
u \circ \left(\bigcup_{i=1}^n U_i\right)
& \supset
uX \cap u \circ \left(uX \cap  \bigcup_{i=1}^n U_i\right)\\
& = {\cls}_\tau \left(uX \cap  \bigcup_{i=1}^n U_i\right)\\
& \supset uX \cap \left(\bigcup_{i=1}^n \pi^{-1}(\pi(U_i))\right).
\end{align*}

Therefore there exist an element $\ga \in \Ga$ and for each $i$ a
set $\{x_{i,j}\}_{j=1}^k \subset U_i$, such that
$d(\ga x_{i,j},x'_{i,j})$ is so small that the inequality
$d(\ga \mu, \mu')< \eta/2$ is satisfied, with
$$
\mu=\frac{1}{nk}\sum_{i=1}^n \sum_{j=1}^k \del_{x_{i,j}}.
$$
Thus $d(\ga\mu,\nu)< \eta$. Since we also have $d(\mu,\ka)< \ep$
and $\ep>0$ is arbitrary,
we have shown that the open set
$$
R_\eta=
\{\mu\in \Mf(X): \ {\text{there exists}}\ \ga\in \Ga\ {\text{ with}}\
d(\ga\mu, \Mf_m(X)) < \eta\}
$$
is dense in $\Mf(X)$.
Clearly $R=\bigcap\{R_\eta: \eta > 0\}$ is the required
dense $G_\del$ subset of $\Mf(X)$.
\end{proof}

\br

\section{Some properties of tame minimal systems}\label{tame-sec}

\begin{thm}\label{tame-inj}(See \cite{G})
Let \XGa be a metric tame dynamical system.
Let $\Mf(X)$ denote the compact convex
set of probability measures on $X$ (with the weak$^*$
topology). Then each element $p\in E(X,\Ga)$ defines
an element $p_*\in E(\Mf(X),\Ga)$ and the map
$p \mapsto p_*$ is both a dynamical system and a semigroup
isomorphism of $E(X,\Ga)$ onto $E(\Mf(X),\Ga)$.
\end{thm}

\begin{proof}
Since $E(X,\Ga)$ is Fr\'echet we have for every $p\in E$ a sequence
$\ga_i\to p$ of elements of $\Ga$ converging to
$p$. Now for every $f\in C(X)$ and every
probability measure $\nu\in \Mf(X)$ we get by the
Riesz representation theorem and
Lebesgue's dominated convergence theorem
$$
\ga_i\nu(f)=\nu(f\circ \ga_i)\to \nu(f\circ p):=p_*\nu(f).
$$
Since the Baire class 1 function $f\circ p$ is well defined
and does not depend upon the choice of the convergent
sequence $\ga_i\to p$, this defines the map $p \mapsto p_*$ uniquely.
It is easy to see that this map is an isomorphism of dynamical systems,
whence a semigroup isomorphism. Finally as $\Ga$ is dense in both
enveloping semigroups, it follows that this isomorphism is onto.
\end{proof}

As we have seen, when \XGa is a metrizable tame system
the enveloping semigroup $E(X,\Ga)$ is a separable Fr\'echet
space. Therefore each element $p\in E$
is a limit of a sequence of elements of $\Ga$,
$p=\lim_{n\to \infty} \ga_n$. It follows that the subset
$C(p)$ of continuity points of each $p\in E$ is a
dense $G_\del$ subset of $X$. More generally, if $Y\subset X$
is any closed subset then the set $C_Y(p)$ of continuity points
of the map $p\rest Y : Y \to X$ is a dense $G_\del$ subset
of $Y$. For an idempotent $v=v^2\in E$ we write
$C_v$ for $C_{\ol{vX}}(v)$.

\begin{lem}\label{Cp}
Let \XGa be a metrizable tame dynamical system,
$E=E(X,\Ga)$ its enveloping semigroup.
\begin{enumerate}
\item
For every $p\in E$ the set $C(p)\subset X$
is a dense $G_\del$ subset of $X$.
\item
For every idempotent $v\in E$, we have
$C_v\subset vX$.
\item
If $v\in E$ is an idempotent such that $\ol{vX}=X$ then
$C(v)\subset vX$.
\item
When $\Ga$ is commutative \XGa minimal we have $C(v)\subset vX$
for every idempotent $v\in E$.
\end{enumerate}
\end{lem}

\begin{proof}
1.\
See the remark above.

2.\
Given $x\in C_v$ choose a sequence $x_n \in vX$ with
$\lim_{n\to\infty} x_n =x$.
We then have $vx=\lim_{n\to\infty} vx_n = \lim_{n\to\infty} x_n =x$,
hence $C_v\subset vX$.

3.\
For such $v$ we have $C(v) = C_v \subset vX$ by part 2.

4.\
When $\Ga$ is Abelian $\ga p =p \ga$ for every
$\ga\in \Ga$ and $p\in E$. In particular the subset
$vX$ is $\Ga$ invariant hence dense in $X$.
\end{proof}

\begin{prop}\label{inv-m}
Let \XGa be a metric tame dynamical system. Then
$p_* \mu =\mu$ for every $p \in E(X,\Ga)$ and every $\Ga$
invariant measure $\mu \in \Mf(X)$.
\end{prop}

\begin{proof}
Let $p = \lim_{n\to\infty} \ga_n$, then for every $f\in C(X)$ and
$x\in X$ we have $\lim_{n\to \infty} f \circ \ga_n (x) = f(px)$, hence
by Lebesgue dominated convergence theorem
\begin{align*}
\int f(x) \,d\mu(x) & =
\lim_{n\to \infty} \int f(\ga_n x) \,d\mu(x)\\
& = \int f(px) \,d\mu(x)\\
&= \int f(x) \,d p_*\mu(x).
\end{align*}
\end{proof}

\begin{prop}\label{pd}
Let $\Ga$ be an Abelian group. Then
any metric tame minimal system \XGa is point distal.
\end{prop}

\begin{proof}
We will prove that the condition in Theorem \ref{vdw} holds; i.e.
that every $S$-set in $X\times X$ is minimal. So let
$W\subset X\times X$ be an $S$-set.
Let $v=v^2$ be some minimal idempotent in $E(X,\Ga)$.
By Theorem \ref{Cp}.3 the set $C(v)$ of continuity points
of the map $v: X \to X$ is a dense $G_\del$ subset of $X$
and moreover $C(v)\subset vX$, so that $vX$ is residual in $X$.
Since by assumption the projection maps $\pi_i: W \to X$
are semiopen, it follows that the sets $\pi_i^{-1}(vX)$
are residual in $W$.
Since $W_{tr}$, the set of transitive points in $W$, is
a dense $G_\del$ subset of $W$ we conclude that
the set $(\pi_1^{-1}(vX) \cap \pi_2^{-1}(vX)) \cap W_{tr}=
(vX \times vX)\cap W_{tr}$ is residual in $W$
and in particular it is nonempty.
Now if $(x,x')$ is any point in this intersection then
$v(x,x') =(vx,vx')=(x,x')$ and $(x,x')$ is a minimal point.
Therefore $W=\OCGa(x,x')$ is minimal.
\end{proof}

\br

\section{Minimal tame systems are almost automorphic and
uniquely ergodic}\label{main-sec}

\begin{thm}\label{almost-auto}
Let $\Ga$ be an Abelian group and \XGa a
metric tame minimal system. Then:
\begin{enumerate}
\item
The system \XGa is almost automorphic. Thus there exist:
\begin{enumerate}
\item
A compact topological group $Y$ with Haar measure $\eta$,
and a group homomorphism
$\kappa:\Ga \to Y$ with dense image.
\item
A homomorphism $\pi:(X,\Ga) \to (Y,\Ga)$,
where the $\Ga$ action on $Y$ is via $\kappa$.
\item
The sets $X_0= \{x\in X: \pi^{-1}(\pi(x)) = \{x\}\}$
and $Y_0=\pi(X_0)$ are dense $G_\del$ subsets of
$X$ and $Y$ respectively.
\end{enumerate}
\item
The system \XGa is uniquely ergodic with unique invariant measure $\mu$
such that $\pi_*(\mu)=\eta$, and $\pi:(X,\mu,\Ga) \to (Y,\eta,\Ga)$is a
measure theoretical isomorphism of the corresponding measure preserving
systems.
\end{enumerate}
\end{thm}

\begin{proof}

1.\
By Proposition \ref{pd} the system $(X,\Ga)$ is point distal.
By the general theory of minimal point distal systems,
if $(X,\Ga)$ is nontrivial,
it admits a unique, nontrivial, maximal equicontinuous factor
$\pi:(X,\Ga) \to (Y,\Ga)$. The system \YGa is thus
of the form stated in the theorem.
If $Y$ is a finite group it follows from structure theory that
$X=Y$. Thus, we now assume that $Y$ is infinite.

Let
\begin{equation*}
\xymatrix
{
X \ar[d]_{\pi}  & X^* =X \vee Y^*\ar[l]_-{\tet^*} \ar[d]^{\pi^*} \\
Y  & Y^*\ar[l]^{\tet}
}
\end{equation*}
be the associated O shadow diagram. Then $(X^*,\Ga)$ is also a minimal
point distal system and \YGa  is its maximal equicontinuous factor.
Again by the theory of point distal systems the extension $\pi^*$ is RIC.
Now, either $\pi^*$ is an isomorphism --- in which case $X=X^*= Y^*$
and $X$ is an almost 1-1 extension of $Y$, hence almost automorphic ---
or there exists a maximal intermediate isometric extension
$$
X^* \overset{\phi}{\to} \hat{X} \overset{\hat\pi}{\to} Y^*
\overset{\tet}{\to} Y,
$$
with $\hat\pi: \hat{X} \to Y^*$ a nontrivial isometric extension.
Thus we now assume that $\hat\pi$ is a nontrivial isometric extension.
As was shown in Proposition \ref{diff} (see also \cite{G1})
$\hat\pi$ has a RIM and we let $y^* \mapsto \la_{y^*}$
be the corresponding section. Associated with this section
we also have
the affine injection $\La: \Mf(Y^*) \to \Mf(\hat X)$ defined by
\begin{equation*}
\La(\nu) =\int_{Y^*} \la_{y^*}\, d\nu(y^*),
\end{equation*}
and we set $\Mf_m(\hat X) :=\{\La(\nu): \nu \in \Mf(Y^*)\}$.

Being isometric, $\hat\pi$ is also open and thus Proposition
\ref{diff} applies (with $\hat X$ and $\hat\pi$ in the roles
of $X$ and $\pi$ respectively).
We conclude that the set $\hat R$ of measures in $\Mf(\hat{X})$ whose orbit
closure meets $\Mf_m(\hat{X}) =\La(\Mf(Y^*))$ is a dense $G_\del$ subset
of $\Mf(\hat{X})$.

\br

2.\
We now have the following diagram
\begin{equation*}
\xymatrix
{
X \ar[dd]_{\pi}  &  X^*\ar[l]_-{\tet^*} \ar[dd]_{\pi^*}
\ar[dr]^{\phi} \\
&  &  \hat X \ar[dl]^{\hat\pi}\\
Y  & Y^*\ar[l]^{\tet}}
\end{equation*}
and we set
\begin{gather*}
\Mf_m(X^*) := \phi^{-1}(\Mf_m(\hat X))
\ {\text{and}}\ \\
\Mf_m(X) :=\tet^*( \Mf_m(X^*))=\tet^*(\phi^{-1}(\Mf_m(\hat X))).
\end{gather*}
Let
\begin{equation*}\label{R*}
R^* :=\phi^{-1}(\hat R) =
\{\xi \in \Mf(\hat X): \text{the orbit closure of $\xi$ meets
$\Mf_m(\hat X)$}\}
\end{equation*}
and
\begin{equation}\label{R}
R' :=\tet^*(R^*) \subset R =
\{\xi \in \Mf(X): \text{the orbit closure of $\xi$ meets
$\Mf_m(X)$}\}.
\end{equation}
By Lemma \ref{so-min} the map $\phi: X^* \to \hat{X}$ is semiopen
and, as $\hat R$ is a dense subset of $\Mf(\hat X)$, Theorem \ref{so-thm}
implies that $R^*$ is a dense subset of $\Mf(X^*)$.
Therefore $R'$ is a dense subset of $\Mf(X)$.
From the definition of $R$ (\ref{R}) it is easy to deduce that it
is a $G_\del$ set and because it contains $R'$,
it is in fact a dense $G_\del$ subset $\Mf(X)$.

\br

3.\
Recall that the system $(X,\Ga)$ is tame
and, by Theorem \ref{tame-inj} so is $(\Mf(X),\Ga)$.
Moreover we have
$E(X,\Ga)=E(\mathfrak{M}(X),\Ga)$.
In particular $u\in E(\mathfrak{M}(X),\Ga)$, as a Baire
class 1 function, has a dense $G_\del$ set of continuity points,
$C_{\mathfrak{M}(X)}(u) \subset \mathfrak{M}(X)$.
Therefore $S:=C_{\mathfrak{M}(X)}(u) \cap R$
is a dense $G_\del$ subset of $\mathfrak{M}(X)$.

Since $uX$ is a dense subset of $X$
it follows that the collection of finite convex combinations
of point masses picked from $uX$ forms a dense
subset of $\Mf(X)$. This implies that
$u\Mf(X)$ is dense in $\Mf(X)$ and
by Lemma \ref{Cp}.3 we have $C_{\mathfrak{M}(X)}(u)
\subset u\mathfrak{M}(X)$.
Thus also $S \subset u\mathfrak{M}(X)$.
Now if $\nu \in S$ then
$u\nu =\nu$ and, $u$ being a minimal idempotent, the
closure of the $\Ga$ orbit of $\nu$ in $\mathfrak{M}(X)$
is a minimal set, whence this entire orbit closure is contained
in $\Mf_m(X)$. In particular $\nu \in \Mf_m(X)$
and we conclude that
$$
S=C_{\mathfrak{M}(X)}(u) \cap R \subset \Mf_m(X).
$$
Therefore $S$ is dense in $\Mf_m(X)$ and in turn, this implies
the equality:
\begin{equation}\label{equ}
\Mf_m(X)=\Mf(X).
\end{equation}

\br

4.\
Given a point $x \in X$, the corresponding point mass
$\del_{x}\in \Mf(X)$ must have, by (\ref{equ}),
a preimage in $\Mf_m(X^*)$, say $\tet^*(\xi)=\del_x$
with $\xi \in \Mf_m(X^*)$.
In particular, for $x$ with ${\tet^*}^{-1}(x)=\{x^*\}$ a singleton,
we must have
$\xi=\del_{x^*} \in \Mf_m(X^*)$ and therefore $\phi_*(\del_{x^*})=
\del_{\hat x}\in \Mf_m(\hat X)$ with $\hat x = \phi(x^*)$.
By the definition of $\Mf_m(\hat X)$ there exists a measure $\rho
\in \Mf(Y^*)$ with
$$
\del_{\hat x} = \La(\rho) =\int_{Y^*} \la_{y^*}\, d\rho(y^*).
$$
This clearly implies that the measures $\rho$ is a point mass,
say $\rho=\del_{y^*}$ and that the measure $\la_{y^*}$ --- which is the
Haar measure on the homogeneous space which forms the fiber
$\hat\pi^{-1}(y^*) \subset \hat X$ --- is also a degenerate point mass.
That is, the isometric extension $\hat\pi$ is in fact an isomorphism.
Now the collapse of $\hat\pi$, implies the collapse of the entire
AI tower, so that in fact $Y^* = X^* =X$, and we have shown that
$X$ is indeed an almost 1-1 extension of $Y$.

\br

5.\
Suppose that $\mu_1$ and $\mu_2$ are two invariant probability
measures on $X$. Then, \XGa being tame, by Proposition \ref{inv-m},
$u_*\mu_i = \mu_i$ and we conclude
that $\mu_i(uX)=1$, for $i=1,2$. Since $\pi$ is
a proximal extension, for every $y\in Y=uY$ the fiber
$\pi^{-1}(y)$ intersects $uX$ at exactly one point:
$\pi^{-1}(y)\cap uX = \{x\}$.
Now by disintegrating each $\mu_i$ over $\eta$, inside the
set $uX$, we conclude that $\mu_1=\mu_2$.
This proves the unique ergodicity of \XGa. It is also clear
from the proof that the map $\pi: (X,\mu,\Ga) \to (Y,\eta,\Ga)$,
where $\mu$ is the unique invariant measure on $X$,
is an isomorphism of measure preserving systems.
\end{proof}

\begin{rmk}
The set $X_0= \{x\in X: \pi^{-1}(\pi(x)) = \{x\}\}$
is a dense $G_\del$ and $\Ga$-invariant subset
of $X$ and thus has $\mu$ measure either zero or one.
In \cite[Section 11]{KL} Kerr and Li construct a minimal Toeplitz
system which is tame and not null.
Since in this construction the growth of the sequence
$\{n_1 < n_2 < \cdots\}$ is arbitrary it follows that the resulting
Toeplitz system can be made not regular in the sense that the densities of
the periodic parts converge to $d < 1$. For such nonregular systems
$\mu(X_0)=0$. This shows that the unique invariant measure of a minimal
tame system need not be supported by the set $X_0$ where $\pi$ is 1-1.
\end{rmk}

\br

%%%%%%%%%%%%%%%%%%%%%%%%%%%%%%%%%%%%%%%%%%%%%%%%%%%%%%%%%%%%%%%%%
%%%%%%             Bibliography                            %%%%%%
%%%%%%%%%%%%%%%%%%%%%%%%%%%%%%%%%%%%%%%%%%%%%%%%%%%%%%%%%%%%%%%%%

\bibliographystyle{amsplain}

\end{document}